\documentclass[11pt]{article}

\usepackage[all]{xy}
\usepackage{graphicx}
\usepackage{amssymb,latexsym,amsmath,amsthm, mathrsfs}
\usepackage{fullpage}

\newcommand{\susp}{\Sigma\mkern-10.1mu/}

\newcommand{\C}{\mathbb C}

\newcommand{\Z}{\mathbb Z}
\newcommand{\N}{\mathbb N}
\def\P{{\mathbb P}}

\def\W{{\mathbb W}}
\def\X{{\mathbb X}}
\def\Y{{\mathbb Y}}
\def\H{{\mathbb H}}
\def\V{{\mathbb V}}

\newtheorem{theorem}{Theorem}[section]
\newtheorem{lemma}[theorem]{Lemma}
\newtheorem{definition}[theorem]{Definition}
\newtheorem{example}[theorem]{Example}

\newtheorem{proposition}[theorem]{Proposition}
\newtheorem{corollary}[theorem]{Corollary}

\newenvironment{acknowledgement}[1][Acknowledgements]
{\begin{trivlist} \item[\hskip \labelsep {\bfseries #1}]}
{\end{trivlist}}

\begin{document}
\title{Semi-topological cycle theory I}

\author{Jyh-Haur Teh}
\date{}
\maketitle \abstract{We study algebraic varieties parametrized by
topological spaces and enlarge the domains of Lawson homology and
morphic cohomology to this category. We prove a Lawson suspension
theorem and splitting theorem. A version of Friedlander-Lawson
moving is obtained to prove a duality theorem between Lawson
homology and morphic for smooth semi-topological projective
varieties. K-groups for semi-topological projective varieties and
Chern classes are also constructed.}

\section{Introduction}
Algebraic cycles are basic ingredients in studying invariants of algebraic varieties. The collection of all $r$-dimensional algebraic cycles of a
projective variety $X$ forms a topological abelian group $Z_r(X)$. Lawson studied these groups from a homotopic viewpoint (\cite{L1}) and proved a
suspension theorem which serves as a cornerstone for Lawson homology and morphic cohomology developed later by Lawson and Friedlander (\cite{F1, FL1,
FL2, FL3}). A continuous map $f:S^n \rightarrow Z_r(X)$ from the $n$-sphere to $Z_r(X)$ can be viewed as a family of algebraic cycles parametrized by
$S^n$. This family can also be considered as ``an" algebraic cycle of $X\times S^n$. This motivates us to consider algebraic varieties parametrized by
topological spaces and consider algebraic cycles on them.

When the base space is an algebraic variety and the parametrization is algebraic, this is just the relative theory of algebraic varieties. The main point
of our studying is that the base space is a very general topological space, the ring of continuous complex-valued functions of it is usually not
Noetherian. It is well known in the algebraic case when we wish a family of algebraic varieties behaves well we need to require the family to be flat
over the base scheme. The flatness of a family of varieties is equivalent to the property that the family is the pull back of the universal family over a
Hilbert scheme by an algebraic morphism to the Hilbert scheme. So to obtain a nice theory, we define our ``semi-topological variety" to be a continuous
map from a topological space $S$ to some Hilbert scheme with some additional technical assumption. We are able to define semi-topological algebraic
cycles on semi-topological projective varieties and extend the definition of Lawson homology and morphic cohomology to them. This paper is the first part
of this theory. A Hodge theory and Riemann-Roch theorem will be given in a forthcoming paper.

Let us give a brief overview of this paper. In section 2, we define semi-topological projective varieties and algebraic cycle on them. Some basic
topological properties of semi-topological cycle groups are studied. In section 3, we prove the Lawson suspension theorem and splitting theorem for
semi-topological projective varieties. In section 4, we give a version of Friedlander-Lawson moving lemma for semi-topological projective varieties and
use it to prove a duality theorem between the Lawson homology and morphic cohomology for semi-topological smooth projective varieties. In section 5, we
compute the Lawson homology group of divisors in a semi-topological smooth projective variety. In section 6, we construct K-groups and Chern classes.

\begin{acknowledgement}
The author thanks Taiwan National Center for Theoretical Sciences(Hsinchu) for proving a nice working environment, and thanks Yu-Wen Kao
for her support.
\end{acknowledgement}

\section{Semi-topological varieties}

Let us briefly recall the construction of cycle groups of complex projective varieties. For a complex projective variety $X$, we write $\mathscr{C}_{r,
d}(X)$ for the collection of effective $r$-cycles of degree $d$ on $X$. According to Chow theorem, $\mathscr{C}_{r, d}(X)$ is a projective variety. Let
$\mathscr{C}_r(X)=\coprod_{d\geq 0}\mathscr{C}_{r, d}(X)$ be the Chow monoid and $Z_r(X)=[\mathscr{C}_r(X)]^+$ the naive group completion of
$\mathscr{C}_r(X)$. Let $K_{r, d}(X)=\coprod_{d_1+d_2\leq d}\pi(\mathscr{C}_{r, d_1}(X)\times \mathscr{C}_{r, d_2}(X))$ where $\pi:\mathscr{C}_r(X)\times
\mathscr{C}_r(X) \rightarrow Z_r(X)$ is the map $(a, b)\mapsto a-b$. We have a filtration
$$K_{r, 0}(X)\subseteq K_{r, 1}(X)\subseteq \cdots =Z_r(X)$$
Each $K_{r, d}(X)$ is compact and the topology of $Z_r(X)$ is the weak topology induced from this filtration. With this topology, $Z_r(X)$ is a
topological abelian group. If $Y$ is also a projective variety, we write $Z_r(Y)(X)$ for the group of algebraic $r$-cocycles on $X$ with values in $Y$,
i.e., $c\in Z_r(Y)(X)$ if $c\in Z_{r+k}(X\times Y)$ where $k$ is the dimension of $X$, the projection from $c$ to $X$ is surjective and fibres of $c$
over $X$ are $r$-cycles in $Y$.

Throughout this paper, $S$ is a compactly generated topological space with based point $s_0$.
We write $\P^n_S$ for $\P^n\times S$.

\begin{definition}
A semi-topological projective variety over $S$ is a continuous map $\X:
S \rightarrow Hil_P(\P^n)$ such that $\X_s$ is a normal projective variety for all $s\in S$ where $Hil_p(\P^n)$
is the Hilbert scheme of $\P^n$ associated to a Hilbert polynomial $p$. We write $\X\subset \P^n_S$ in this case. We define the dimension of $\X$ to be the dimension of $\X_{s_0}$, and $|\X_{s_0}|$ for the algebraic variety corresponding to $\X_{s_0}$.
\end{definition}

\begin{definition}
Suppose $\X\subset \P^n_S, \Y\subset \P^m_S$ are semi-topological projective varieties over $S$.
Let
$$Z_r(\Y)(\X):=\{\alpha\in Map((S, s_0), (Z_k(\P^n\times \P^m), 0))|\alpha(s)\in Z_r(\Y_s)(\X_s)\}$$ where
$k=r+dim\X$. Let $Z_{k, \leq e}(\P^n\times \P^m)$ be the image of
$\coprod_{e_1+e_2\leq e}\mathscr{C}_{k, e_1}(\P^n\times \P^m)\times \mathscr{C}_{k, e_2}(\P^n\times \P^m)$
in $Z_k(\P^n\times \P^m)$. The topology of $Z_k(\P^n\times \P^m)$ is the weak topology induced from the filtration
$$Z_{k, \leq 0}(\P^n\times \P^m) \subseteq Z_{k, \leq 1}(\P^n\times \P^m) \subseteq \cdots =Z_k(\P^n\times \P^m)$$
Let $Z_{r, \leq e}(\Y)(\X):=Map((S, s_0), (Z_{k, \leq e}(\P^n\times \P^m), 0))\cap Z_r(\Y)(\X)$. We give
$Z_r(\Y)(\X)$ the weak topological induced from the filtration
$$Z_{r, \leq 0}(\Y)(\X)\subseteq Z_{r, \leq 1}(\Y)(\X) \subseteq \cdots =Z_r(\Y)(\X)$$
where each $Z_{r, \leq e}(\Y)(\X)$ is endowed with the compact-open topology. Then
$Z_r(\Y)(\X)$ is a topological abelian group.
\end{definition}

\begin{definition}
If $\Y, \Y'$ are two semi-topological projective varieties, we say that $\Y'$ is a subvariety of $\Y$,
denoted as $\Y'\subseteq \Y$
if $|\Y'_s| \subseteq |\Y_s|$ for all $s\in S$. A semi-topological Zariski open set if a set of the form $|\Y-\Y'|$ where $\Y'\subseteq \Y$.
\end{definition}

\begin{definition}
Suppose $\X\subset \P^n_S$, $\Y\subset \P^m_S$. Suppose for each $s\in S$, $f_s:|\X_s| \rightarrow |\Y_s|$ is a
given morphism of projective varieties. The assignment $s \overset{f}{\mapsto} f_s$ is said to be a
morphism between $\X$ and $\Y$,
denoted $f:\X \rightarrow \Y$ if $s\mapsto gr f_s\in \mathscr{C}_r(\P^n\times \P^n)$ is continuous
where $gr f_s$ is the graph of $f_s$.
\end{definition}

\begin{definition}
If $f:\X \rightarrow \Y$ is a morphism of semi-topological projective varieties, define
$f_*:Z_r(\X)(\W) \rightarrow Z_r(\Y)(\W)$ by
$$(f_*\alpha)(s):=q_{s*}((|\W_s|\times gr f_s)\bullet p^*_s(\alpha(s)))$$
where $p:\W\times \X\times \Y \rightarrow \W\times \X, q:\W\times \X\times \Y \rightarrow \W\times \Y$
are projections.
\end{definition}

\begin{proposition}
$f_*$ is continuous.
\end{proposition}

\begin{proof}
We note that it is well known that the intersection product on cycles intersecting properly is continuous
\cite{Fu}.
Suppose $\beta:T \rightarrow Z_r(\X)(\W)$ is continuous where $T$ is some topological space. Then
$(f_*(\beta(t)))(s)=q_{s*}((|\W_s|\times gr f_s)\bullet p^*_s(\beta(t)(s)))$ which is continuous in $t$.
Thus $f_*$ is continuous.
\end{proof}

\begin{proposition}
If $\Y'$ is a subvariety of $\Y$, then $Z_r(\Y')(\X)$ is closed in $Z_r(\Y)(\X)$.
\end{proposition}

\begin{proof}
For $\alpha\in Z_{r, \leq e}(\Y)(\X)-Z_{r, \leq e}(\Y')(\X)$, there
is $s_1\in S$ such that $\alpha(s_1)\in Z_{r, \leq
e}(\Y'_{s_1})(\X_{s_1})-Z_{r, \leq e}(\Y'_{s_1})(\X_{s_1})$. Since
$Z_{r, \leq e}(\Y'_{s_1})(\X_{s_1})$ is closed in $Z_{r, \leq
e}(\Y_{s_1})(\X_{s_1})$, see \cite[Prop 2.9]{T3}, there is $V$ open in $Z_{k,
\leq e}(\P^n\times \P^m)$ such that $\alpha(s_1)\in V\cap Z_{r, \leq
e}(\Y_{s_1})(\X_{s_1})\subset Z_{r, \leq
e}(\Y_{s_1})(\X_{s_1})-Z_{r, \leq e}(\Y'_{s_1})(\X_{s_1})$. Let
$$W=\{\beta\in Map((S, s_0), (Z_{k, \leq e}(\P^n\times \P^m),
0))|\beta(s_1)\in V\}$$ then $W$ is open in $Map((S, s_0), (Z_{k,
\leq}(\P^n\times \P^m), 0))$ and $\alpha\in W\cap Z_{r, \leq
e}(\Y)(\X)\subseteq Z_{r, \leq e}(\Y)(\X)-Z_{r, \leq e}(\Y')(\X)$.
Therefore $Z_{r, \leq e}(\Y)(\X)-Z_{r, \leq e}(\Y')(\X)$ is open and
$Z_{r, \leq e}(\Y')(\X)$ is closed in $Z_{r, \leq e}(\Y)(\X)$. We
then have $Z_r(\Y')(\X)$ is closed in $Z_r(\Y)(\X)$.
\end{proof}

Recall that there is a functor $k$ (see \cite{S}) constructed by
Steenrod from the category of topological spaces to the category of
compactly generated spaces which acts like a retraction.
Furthermore, for any topological space $X$, $X$ and $k(X)$ have the
same homology and homotopy groups. Recall that by the construction
in \cite{T2}, if $H$ is a normal closed subgroup of $G$ and both are
compactly generated, then the short exact sequence $0\rightarrow H
\rightarrow G \rightarrow G/H \rightarrow 0$ gives a fibration
$$\xymatrix{B_H \ar[r] & B_G \ar[d]\\
            & B_{G/H}}$$
where $BG$ is the classifying space of $G$.
Thus we have a long exact sequence of homotopy groups
$$\cdots \rightarrow \pi_n(H) \rightarrow \pi_n(G) \rightarrow
\pi_n(G/H) \rightarrow \pi_{n-1}(H) \rightarrow \cdots$$ Combine
with the Steenrod functor $k$, once we have some complicated
topological abelian groups that form the short exact sequence stated
above, we get a long exact sequence of homotopy groups. The
following is an application of this result.

\begin{definition}
Let $$Z_r(\Y; \Y')(\X):=\frac{Z_r(\Y)(\X)}{Z_r(\Y')(\X)}$$ and
$$Z^t(\X):=Z_r(\P^t_S; \P^{t-1}_S)(\X)$$ where $\Y'$ is a semi-topological subvariety
of $\Y$.
\end{definition}

\begin{corollary}
We have a long exact sequence of homotopy groups
$$\cdots \rightarrow \pi_nZ_r(\Y')(\X) \rightarrow \pi_nZ_r(\Y)(\X)
\rightarrow \pi_nZ_r(\Y; \Y')(\X) \rightarrow \pi_{n-1}Z_r(\Y')(\X) \rightarrow \cdots$$
\end{corollary}

\begin{definition}
Let $pt:S \rightarrow Hil_p(\P^n)$ be a constant map whose image is a
point in $\P^n$. Then $Z_r(\Y)(pt)$ is isomorphic to $Z_r(\Y)(pt')$
for any two such maps $pt, pt'$. We write $Z_r(\Y):=Z_r(\Y)(pt)$
without referring to which point we take. The map $pt$ is called a
point map.
\end{definition}

\begin{definition}
Define $$\H_{S, n}(\X):=\pi_nZ_0(\X)$$ and
$$\H^n_S(\X):=\pi_{2m-n}Z^m(\X)$$ where $m$ is the dimension of $\X$.
\end{definition}

\begin{example}
When $S=S^0$ the $0$-dimensional sphere, $\X=X\times S^0$ for some projective variety $X$,
$$\H_{S, n}(\X)=\pi_nZ_0(X)\cong H_n(X)$$ by the Dold-Thom theorem. If $X$ is smooth, then
$$\H^n_S(\X)=\pi_{2m-n}Z^m(X)\overset{FL}{\cong} \pi_{2m-n}Z_0(X)\cong H_{2m-n}(X)\overset{PD}{\cong} H^n(X)$$
where $FL$ is the Friedlander-Lawson duality isomorphism and $PD$ is the Poincare duality isomorphism.
\end{example}

\section{Suspension theorem and splitting theorem}
Let us recall that if $X\subseteq \P^n$ and $x_{\infty}\in \P^{n+1}\backslash \P^n$, the suspension $\susp X$ of $X$ with respect to $x_{\infty}$ is the join of $X$ and $x_{\infty}$.

\begin{definition}
Suppose that $\Y\subseteq \P^m_S$, $pt$ a point map with image in
$\P^{m+1}\backslash \P^m$. The suspension of $\Y$ with respect to
$pt$ is the semi-topological subvariety
$$(\susp_{pt}\Y)(s):=\susp_{pt(s)}\Y_s$$
\end{definition}
So we have $\susp_{pt}\Y\subseteq \P^{m+1}_S$. The suspension
induces a map $Z_k(\P^n\times \P^m)
\overset{\susp_{pt}}{\longrightarrow} Z_{k+1}(\P^n\times \P^{m+1})$
by suspending $\P^m$. Hence for each $\alpha\in Z_r(\Y)(\X)$,
$\susp_{pt}$ induces a semi-topological cycle in
$Z_{r+1}(\susp\Y)(\X)$ in the following way
$$(\susp_{pt}\alpha)(s):=\susp_{pt}\alpha(s)$$

\begin{theorem}
Let $pt=[0:\cdots:0:1]\in \P^{m+1}, \Y\subseteq \P^m_S, \X\subseteq
\P^n_S$. Then $\susp_{pt*}:Z_r(\Y)(\X) \rightarrow
Z_{r+1}(\susp_{pt}\Y)(\X)$ is a weak homotopy equivalent.
\end{theorem}

\begin{proof}
We write $\susp$ for $\susp_{pt}$. Let
$$T_{r+1}(\susp \Y)(\X):=\{\alpha \in Z_{r+1}(\susp \Y)(\X)|\alpha(s) \mbox{ meets }
\X_s\times \Y_s \mbox{ properly in } \X_s\times \susp\Y_s, \mbox{
for all }s\in S\}$$ Let $\Lambda\subseteq \P^{m+1}\times \P^1\times
\P^{m+1}$ be the closed subvariety constructed by Friedlander in
\cite[Prop 3.2]{F1} which is a geometric description of Lawson's
holomorphic taffy. Let $\Lambda_t:=\Lambda\bullet (\P^{m+1}\times
\{t\}\times \P^{m+1}), t\in \C$. Then for $\alpha \in T_{r+1}(\susp
\Y)(\X)$,
$$\Phi_t(\alpha):=q_{t*}(p^*_t\alpha\bullet (\P^n\times
\Lambda_t))\in T_{r+1}(\susp \Y)(\X)$$ where $p_t:\P^n\times
\P^{m+1}\times \{t\}\times \P^{m+1} \rightarrow \P^n\times \P^{m+1}$
and $q_t:\P^n\times \P^{m+1}\times \{t\}\times \P^{m+1} \rightarrow
\P^n\times \P^{m+1}$ are the projections to the $(1, 2)$-component
and the $(1, 4)$-component respectively. If $t=0$,
$\Phi_0(\alpha)\in \susp Z_r(\Y)(\X)$. Not difficult to see that
$\Phi$ is a strong deformation retract of $T_{r+1}(\susp \Y)(\X)$ to
$\susp Z_r(\Y)(\X)$.

Let $x_1=[0:\cdots :0:1]\in \P^{m+2}, x_2=[0:\cdots:0:1:1]\in
\P^{m+2}$. Recall that by \cite[Prop 2.3]{FL2}, for any $d>0$, there is $e(d)>0$
such that for any $e>e(d)$, there is a line $L_e$ in $\mathscr{C}_{m+1,
e}(\P^{m+2})$ containing $e\P^{m+1}$ such that we have a map
$$\Psi_e:Z_{r+1, \leq d}(\P^{m+1})\times L_e \rightarrow Z_{r+1,
\leq de}(\P^{m+1})$$ defined by
$$\Psi_e(Z, D):=p_{2*}((x_1\#Z)\bullet D)$$
where $p_2:\P^{m+2}-\{x_2\} \rightarrow \P^{m+1}$ is the projection
with center $\{x_2\}$. Furthermore, for $D\in L_e-\{e\P^m\},
\Psi_e(Z, D)\in T_{r+1, de}(\P^{m+1})$ and $\Psi_e(Z, e\P^m)=eZ$.
When we restrict cycles having support in $\susp \Y\subseteq
\P^{m+1}$ by checking the definition of $\Psi_e$, we get a map
$$\Psi_e:Z_{r+1, \leq d}(\susp \Y)\times L_e \rightarrow Z_{r+1,
\leq de}(\susp \Y)$$ with the corresponding properties. For $\alpha\in Z_{r+1, \leq d}(\susp \Y)(\X)$, define
$$\Psi_e(Z, D)(s):=p_{2*}((x_1\#\alpha(s))\bullet (D\times |\X_s|))$$
This map is continuous in $s$ and have the homotopy property as
before. Note that if $f:C \rightarrow Z_{r+1}(\susp \Y)(\X)$ is a
map from a compact topological space $C$, the image $Im f$ is
compact and $Im f\subseteq Z_k(\P^n\times \P^m)$, by \cite[Lemma 2.8]{T3}, $Im
f\subseteq Z_{k, \leq d}(\P^n\times \P^m)$ for some $d>0$. Therefore
$Im f\subseteq Z_{r+1, \leq d}(\susp \Y)(\X)$.

We show that the map
$i_*:T_{r+1}(\susp \Y)(\X) \rightarrow Z_{r+1}(\susp \Y)(\X)$
induced from the inclusion is a weak homotopy equivalence.
Let $[f]\in \pi_n(Z_{r+1}(\susp \Y)(\X))$ be a based point preserving continuous map. Since $Im f$ is compact, $Im f\subseteq Z_{r+1, \leq d}(\susp \Y)(\X)$ for some $d$. Then by the result above, there is a map $\Psi_e:Z_{r+1, \leq d}(\susp \Y)(\X)\times L_e \rightarrow Z_{r+1, \leq de}(\susp \Y)(\X)$ such that $\Psi_e(f(s), e\P^m)=f(s)$, $\Psi_e(f(s), D)\in T_{r+1}(\susp \Y)(\X)$ for $D\in L_e-\{e\P^m\}$. Hence $i_*[\Psi_e(f, D)]=[f]$ which implies that $i_*$ is surjective. For injectivity, if $[g]\in \pi_nT_{r+1}(\susp \Y)(\X)$ is mapped to 0 by $i_*$ in $Z_{r+1}(\susp \Y)(\X)$, then $i_*g$ can be extended to a map $\widetilde{g}:D^{n+1} \rightarrow Z_{r+1}(\susp \Y)(\X)$ where $D^{n+1}$ is the unit closed ball. Again, by choose some $\Psi_e$, we can show that $\widetilde{g}$ is homotopy to some $\Psi_e(\widetilde{g}, D):D^{n+1} \rightarrow T_{r+1}(\susp \Y)(\X)$. Thus $[g]=0$. Combining
with previous result, the proof is complete.
\end{proof}

\begin{theorem}(The splitting theorem)
If $\X$ is a semi-topological projective variety, then there is
$$\xi_t:Z_0(\P^t)(\X) \rightarrow Z^t(\X)\times Z^{t-1}(\X)\times
\cdots \times Z^0(\X)$$ which is a weak homotopy equivalence.
\end{theorem}

\begin{proof}
Recall that there is an isomorphism $\P^n\cong \mathscr{C}_{0,
n}(\P^1)$ for any positive integer $n$. The projection map
$\P^t\cong \mathscr{C}_{0,t}(\P^1) \rightarrow \mathscr{C}_{0,
\binom{t}{k}}(\P^k)$ defined by
$$x_1+\cdots +x_t\mapsto \sum_{I \subset \{1, ..., t\}, |I|=k}x_I$$
where $x_I=x_{i_1}+\cdots +x_{i_t}$ for $I=\{i_1, ..., i_t\}$ induces a
map $\xi_t:Z_0(\P^t_S)(\X) \rightarrow Z_0(\P^k_S)(\X) \rightarrow
Z^k(\X)$ for $0\leq k\leq t$. We have a commutative diagram
$$\xymatrix{Z_0(\P^{t-1}_S)(\X) \ar[d] \ar[r]^{\xi_{t-1}} &
Z^{t-1}(\X)\times \cdots \times Z^0(\X) \ar[d] \\
Z_0(\P^t_S)(\X) \ar[d]_{q} \ar[r]^-{\xi_t} & Z^t(\X)\times Z^{t-1}(\X) \times \cdots \times Z^0(\X)\ar[d]^{p}\\
Z^t(\X) \ar[r]^{=} & Z^t(\X)}$$ where $q$ is the quotient map, $p$
is the projection map. From the homotopy sequence associated to the
vertical columns, we get the result by induction on $t$.

\end{proof}

\section{Moving lemma}
Let $H=Hil_p(\P^n)$ be the Hilbert scheme of $\P^n$ associated to the Hilbert polynomial $p$, and $\widetilde{H}\overset{\pi}{\longrightarrow} H$
be the universal family over $H$. Suppose each algebraic variety parametrized by $H$ is of dimension $m$.
Let $U_{\widetilde{H}}(d)\subset \P(\Gamma(\mathcal{O}_{\P^n_H}(d)^{m+1}))$ be the Zariski open set of those
$F=(f_0, ..., f_m)$ such that
$$\mathscr{L}_F:=\{(t, h)\in \P^n_H|F_h(t)=0\}$$ misses $\widetilde{H}$ where $F_h=(f_{0, h}, ..., f_{m, h})$ is obtained from pulling back
$F$ by the inclusion $z\mapsto (z, h)$ from $\P^n$ to $\P^n_H$. Then $F$ induces a finite morphism $p_F:\widetilde{H} \rightarrow \P^m_H$
by $p_F(x):=p_{F_{\pi(x)}}(x)$.

For $Y\in \mathscr{C}_{r, \leq e}(\widetilde{H})(H), Z\in \mathscr{C}_{\ell, \leq e}(\widetilde{H})(H)$, let
$$Y\star_F Z:=\{(y, z)\in Y\times_H Z|y\neq z, p_F(y)\neq p_F(z)\}$$
where $Y\times_H Z$ is the fibre product of $Y$ and $Z$ over $H$.

Follow similar approach as in \cite[Prop 1.3]{FL2}, we get the following result
\begin{proposition}
Suppose that $r+\ell\geq m$, $e\in \N$. There is a Zariski closed subset $\mathscr{B}(d)_e\subset U_{\widetilde{H}}(d)$ with
$\lim_{d \to \infty}Fcodim\mathscr{B}(d)_e=\infty$ where $Fcodim\mathscr{B}(d)_e=\underset{h}{min}\{codim \mathscr{B}(d)_{e, h}\}$ and
$\mathscr{B}(d)_{e, h}:=\{F_h|F\in \mathscr{B}(d)_{e, h}\}$ such that for any $Y\in \mathscr{C}_{r, \leq e}(\widetilde{H})(H),
Z\in \mathscr{C}_{\ell, \leq e}(\widetilde{H})(H)$, $|Y_h|\star_{F_h}|Z_h|$ has pure dimension $r+\ell-m$ whenever $F\in U_{\widetilde{H}}(d)-\mathscr{B}(d)_e$.
\end{proposition}

Now let $\X:S \rightarrow H$ be a semi-topological projective variety. Then the pullback $\X^*\P^n_H=\P^n_S$ and for $F\in \P(\Gamma(\mathcal{O}_{\P^n_H}(d))^{m+1})$,
$\X^*F(x, s):=F(x, \X(s)), \X^*F\in \P(\Gamma(\mathcal{O}_{\P^n_S}(d)^{m+1}))$. Define
$p_{\X^*F}:\X^*(\widetilde{H}) \rightarrow \P^m_H$ by
$$p_{\X^*F}(x, s):=p_{F_{\X(s)}}(x)$$ where $x\in \widetilde{H}, s\in S$.
For $\alpha\in Z_{r, \leq e}(\X), \beta\in Z_{\ell, \leq e}(\X)$,
$$\alpha\star_{\X^*F}\beta:=\{(a, b)\in |\alpha|\times_S|\beta||a\neq b, p_{\X^*F}(a, s)\neq p_{\X^*F}(b, s)\}$$
Let $U_{\widetilde{H}}(d)\subset \P(\Gamma(\mathcal{O}_{\P^n_S}(d)^{m+1}))$ be the semi-topological Zariski open set of those
$F=(f_0, ..., f_m)$ such that
$$\mathscr{L}_F:=\{(t, s)\in \P^n_S|F_s(t)=0\}$$ misses $\X^*\widetilde{H}$.
Then by taking $\widetilde{\mathscr{B}}(d)_e=\X^*(\mathscr{B}(d)_e)$, from the above result, we also have enough good projections for semi-topological projective varieties when the degree is large enough.

\begin{corollary}
Let $\X\subseteq \P^n_S$ be a semi-topological projective variety of dimension $m$.
Suppose that $r+\ell\geq m$, $e\in \N$. There is a semi-topological Zariski closed subset $\widetilde{\mathscr{B}}(d)_e\subseteq U_{\X}(d)$
with $\lim_{d\to \infty}Fcodim\mathscr{B}(d)_e =\infty$ such that
for $\alpha\in Z_{r, \leq e}(\X), \beta\in Z_{\ell, \leq e}(\X)$, $|\alpha_s|\star_{F_s}|\beta_s|$ has pure dimension $r+\ell-m$ where
$F\in U_{\X}(d)-\widetilde{\mathscr{B}}(d)_e$.
\end{corollary}

Once we know how to find good projections for semi-topological projective varieties, follow argument of Friedlander and Lawson in \cite{FL2}, we get a moving lemma for semi-topological projective varieties.

\begin{theorem}
Let $\X\subseteq \P^n_S$ be a semi-topological projective variety of dimension $m$. Let $r, \ell, e$ be nonnegative integers with $r+\ell\geq m$. Then
there exist an open set $\mathcal{O}$ of $\{0\}$ in $C$ and a continuous map
$$\widetilde{\Psi}:\mathscr{C}_{\ell}(\X)\times \mathcal{O} \rightarrow \mathscr{C}_{\ell}(\X)^2$$ such that
$\pi\circ \widetilde{\Psi}$ induces by linearity a continuous map
$$\Psi:Z_{\ell}(\X)\times \mathcal{O} \rightarrow Z_s(\X)$$ satisfying the following properties. Let $\psi_p=\Psi|_{Z_s(\X)\times \{p\}}$ for $p\in \mathcal{O}$.
\begin{enumerate}
\item $\psi_0=Id$.
\item For any $p\in \mathcal{O}$, $\psi_p$ is a continuous group homomorphism.
\item For any $\alpha\in Z_{\ell, \leq e}(\X)$, $\beta\in Z_{r, \leq e}(\X)$, and any $p\neq 0$ in $\mathcal{O}$, each component of excess dimension
of the intersection $|\alpha(s)|\cap |\psi_p(Z)|$ is contained in the singular locus of $|\X_s|$, for $s\in S$.
\end{enumerate}
\end{theorem}

Let $(\X, \Y)$ be a pair of semi-topological projective varieties in $\P^n_S$ where $\Y\subseteq \X$. We say that
a map $f:(\X, \Y) \rightarrow (\X', \Y')$ between two pairs of semi-topological varieties is a relative isomorphism
if $f:\X \rightarrow \X'$ is a semi-topological morphism such that $f:\X-\Y \rightarrow \X'-\Y'$ is an isomorphism of semi-topological quasi-projective varieties. The following example is the most important case to us.

\begin{example}
Define $\phi:(\X\times_S \P^t_S, \X\times_S \P^{t-1}_S) \rightarrow (\susp^t\X, \P^{t-1}_S)$ by
$$\phi(([x_0:\cdots :x_n], s), ([a_0:\cdots:a_t], s)):=([a_0x_0:\cdots:a_0x_n:a_1:\cdots:a_t], s)$$ where we identify the $\P^{t-1}_S$ of the second pair
to the hyperplane at infinity of $\susp^t\X$. Then not difficult to see that $\phi$ is a relative isomorphism.
\end{example}

The following lemma is a special case in order to define the cycle groups for quasi-projective variety (see \cite{Li1}) but it is enough to prove the duality theorem.

\begin{lemma}
Suppose that $f:(\X, \Y) \rightarrow (\X', \Y')$ is a relative isomorphism where $dim\Y'<r$, then
$\frac{Z_r(\X)}{Z_r(\Y)}$ is weak homotopic equivalent to $\frac{Z_r(\X')}{Z_r(\Y')}$ and $\frac{Z_r(\X')}{Z_r(\Y')}$ is isomorphic
to $Z_r(\X')$.
\end{lemma}

\begin{proof}
The morphism $f:\X \rightarrow \X'$ induces group homomorphisms $f_*:Z_r(\X) \rightarrow Z_r(\X')$ and $f_*:Z_r(\Y) \rightarrow Z_r(\Y')$.
Since $f$ restricts to $\X-\Y$ is injective, this gives the injectivity of $f_*:\frac{Z_r(\X)}{Z_r(\Y)} \rightarrow \frac{Z_r(\X')}{Z_r(\Y')}$. But
since $r>dim \Y'$, $Z_r(\Y')=\{0\}$. Hence $f_*$ is surjective and $\frac{Z_r(\X')}{Z_r(\Y')}=Z_r(\X')$.
\end{proof}

By using the moving lemma, we got the following duality theorem which is proved by similar arguments in \cite{FL3}.

\begin{theorem}(Duality theorem)
Suppose that $\X\subseteq \P^n_S, \Y\subseteq \P^k_S$ where $dim\X=m$, $\X, \Y$ are smooth. Then there is a weak homotopic equivalence:
$i_*:\Z_k(\Y)(\X) \cong \Z_{m+k}(\X\times \Y)$ where $i$ is the inclusion.
\end{theorem}

\begin{corollary}
Suppose that $\X$ ia a nonsingular semi-topological projective variety with dimension $m$. If $0\leq t \leq m$,
then $Z^t(\X)$ is weak homotopic equivalent to $Z_{m-t}(\X)$.
\end{corollary}

\begin{proof}
$Z^t(\X)=\frac{Z_0(\P^t)(\X)}{Z_0(\P^{t-1})(\X)}\cong \frac{Z_m(\X\times \P^t_S)}{Z_m(\X\times \P^{t-1}_S)}\cong \frac{Z_m(\susp^t\X)}{Z_m(\X\times \P^{t-1})}\cong Z_m(\susp^t\X)\cong Z_{m-t}(\X)$
\end{proof}

\section{Semi-topological divisors}
Suppose that the dimension of $\X$ is greater than 0. Let
$$K[\X]=\mathscr{C}(S)[Z_0, ...., z_n]/I(\X)=\bigoplus_{d\geq
0}K_d(\X)$$ where $\mathscr{C}(S)$ is the ring of complex-valued
continuous functions on $S$, and $K_d(\X)$ is the collection of
homogeneous polynomials of degree $d$ in $K[\X]$.

\begin{proposition}
If $f\in \mathscr{C}(S)[z_0, ..., z_n]$ is homogenous of degree $d$
which is not a zero polynomial for any $s\in S$, then $f$ defines an effective
semi-topological divisor $(f)\in \mathscr{C}_{n-1}(\P^n_S)$ by
$$(f)(s):=(f_s)$$ for $s\in S$ where $(f_s)$ is the divisor on
$\P^n$ defined by $f_s$.
\end{proposition}

\begin{proof}
Since $f_s$ is not a zero polynomial for any $s\in S$, $(f_s)$ is an
effective divisor for any $s\in S$. From the definition of Chow
form, we see that the coefficients of the Chow form $F_{(f_s)}$ of
$(f_s)$ are continuous functions of the coefficients of $f$. This
implies that the assignment $(f):S \rightarrow \mathscr{C}_{n-1,
d}(\P^n)$ is continuous.
\end{proof}

\begin{definition}
Let $\mathscr{C}(S)[z_0, ..., z_n]_{d, \X}$ be the collection of all
$f\in \mathscr{C}(S)[z_0, ..., z_n]$ of degree $d$ such that $(f_s)$ meets $\X_s$ properly in $\P^n$ for all $s\in S$.
For $f+I(\X)\in K_d(\X)$ where $f\in \mathscr{C}(S)[z_0, ..., z_n]_{d, \X}$, let
$$(f+I(\X))(s):=(f_s)\bullet \X_s$$ for $s\in S$
Then $(f+I(\X))$ is a semi-topological divisor on $\X$. Let
$$W_d(\X):=\{(f+I(\X))|f\in \mathscr{C}(S)[z_0, ..., z_n]_{d, \X}\}$$
Let $W(\X)=\coprod_{d\geq 0}W_d(\X)$ and
$$Z_{m-1}(\X)^{lin}=\{\alpha-\beta|\mbox{ where } \alpha, \beta\in W(\X), \alpha(s_0)=\beta(s_0)\}$$  We say that
a semi-topological divisor $D\in Z_{m-1}(\X)$ is semi-topologically linearly equivalent to zero if
$D\in Z_{m-1}(\X)^{lin}$.
\end{definition}

\begin{proposition}
Let $T_d(\X)=\{(f)|f\in \mathscr{C}(S)[z_0, ..., z_n]_{d, \X}\}$,
$T(\X)=\coprod_{d\geq0}T_d(\X)$ and
$$\widetilde{T}(\X):=\{(f)-(g)|(f)(s_0)=(g)(s_0), (f), (g)\in T(\X)\}\subseteq Z_{m-1}(\X)$$

 Then
\begin{enumerate}
\item $\widetilde{T}(\X)$ is isomorphic as a topological group to $Z_{m-1}(\X)^{lin}$.

\item
$\widetilde{T}(\X)$ is weak homotopy equivalent to
$Z_{n-1}(\P^n_S)$ where $\X\subseteq \P^n_S$.
\end{enumerate}
\end{proposition}

\begin{proof}
The isomorphism between $\widetilde{T}(\X)$ and $Z_{m-1}(\X)^{lin}$ is given by the natural map$(f)-(g)\mapsto (f+I(\X))-(g+I(\X))$.
By moving lemma, for any $e>0$, there is an integer
$e(d)$ such that for any $k>e(d)$ there is a continuous function
$\Theta_{k}:Z_{n-1, \leq e}(\P^n)\times \ell^0 \rightarrow
Z_{n-1, ke}(\P^n)$ such that
\begin{enumerate}
\item $\Theta_{k}(c, 0)=kc$,
\item $\Theta_{k}(c, t)$ meets $\X_s$ properly for $t\in\ell^0\backslash \{0\}$.
\end{enumerate}

Then follow exactly the argument in proving the suspension theorem, we show that the inclusion
$i_*:\widetilde{T}(\X) \rightarrow Z_{n-1}(\P^n_S)$ is a weak homotopy equivalence.
\end{proof}

\begin{proposition}
Suppose that $\X\subseteq \P^n_S$ is a semi-topological variety of
dimension $m$. Then $Z_{m-1}(\X)^{lin}$ is weak homotopy equivalent
to $Map((S, s_), (Z_0(\P^1), 0))$. In particular,
$$\pi_{\ell}Z_{m-1}(\X)^{lin}=\left\{
                           \begin{array}{ll}
                             H^2(S), & \hbox{ if }\ell=0 \\
                             H^1(S), & \hbox{ if } \ell=1\\
                             H^0(S), & \hbox{ if }\ell=2\\
                             0, & \hbox{ otherwise}
                           \end{array}
                         \right.$$
\end{proposition}

\begin{proof}
We have a homeomorphism between $\mathscr{C}(S)[z_0, ..., z_n]_d$ and
$Map(S, \C^{\binom{n+d}{d}})$
by $$f\mapsto \mbox{ coefficients of } f$$ This homeomorphism reduces to a homeomorphism
$$\mathscr{C}_{n-1, d}(\P^n) \cong Map(S, \P^{\binom{n+d}{d}-1})\cong Map(S, \mathscr{C}_{0, \binom{n+d}{d}-1}(\P^1))$$

Thus $Z_{n-1}(\P^n_S)\cong Map((S, s_0), (Z_0(\P^1), 0))$. From the result above, we have weak homotopic equivalent $Z_{m-1}(\X)^{lin}\cong \widetilde{T}(\X)\cong Z_{n-1}(\P^n_S)\cong map((S, s_0), (Z_0(\P^1), 0)$.
\end{proof}

\section{Chern classes}
\begin{definition}
Suppose that the dimension of a semi-topological variety $\X$ is
$k$. Let
$$\mathscr{C}^{s, 1}(\P^n_S)(\X):=\{\alpha\in Map(S,
\mathscr{C}_{k+n-s}(\P^m\times \P^n))|\alpha(s)\in \mathscr{C}^{s,
1}(\P^n)(|\X_s|)\}$$ where $\mathscr{C}^{s, 1}(\P^n)(|\X_s|):=\mathscr{C}_{n-s, 1}(|\X_s|\times \P^n)$.
\end{definition}

By suspension, we have a sequence
$$\cdots \rightarrow \mathscr{C}^{s, 1}(\P^n)(\X) \rightarrow \mathscr{C}^{s,
1}(\P^{n+1})(\X) \rightarrow \mathscr{C}^{s,
1}(\P^{n+2})(\X)\rightarrow \cdots$$ Let
$$\mathscr{C}^{s, 1}(\P^{\infty})(\X):=\lim_{n\to \infty}\mathscr{C}^{s,
1}(\P^n)(\X)$$ and let
$$\V ect^s(\X):=[\mathscr{C}^{s, 1}(\P^{\infty})(\X)]^+$$ be the
group completion. Note that we do not fix a based point for $\V
ect^s(\X)$. Let
$$\widetilde{\V ec}t^s(\X):=\{f-g\in \V ect^s(\X)|f_{s_0}=g_{s_0}\}$$
and let
$$\widetilde{\V ec}t^s(\X)_n:=\{f-g\in \widetilde{\V}ect^s(\X)|f,
g\in \mathscr{C}^{s, 1}(\P^n)(\X)\}$$ then we have the following sequences and maps
$$\xymatrix{ \widetilde{\V ec}t^s(\X)_n \ar[d]_{\susp} \ar[r] &
Z_{n-s}(\P^n_S)(\X) \ar[d]^{\susp}\\
\widetilde{\V ec}t^s(\X)_{n+1} \ar[d]_{\susp} \ar[r] &
Z_{n-s}(\P^{n+1}_S)(\X) \ar[d]^{\susp}\\
\vdots & \vdots}$$ and we get a map
$$\widetilde{\V ec}t^s(\X) \rightarrow \lim_{n\to
\infty}Z_{n-s}(\P^n_S)(\X)$$ By taking $\pi_0$ on both sides, we get a homomorphism
$$\pi_0\widetilde{\V ec}t^s(\X) \overset{\textbf{c}}{\longrightarrow}
\pi_0(\lim_{n\to \infty}Z_{n-s}(\P^n_S)(\X))\cong \pi_0Z_o(\P^n_S)(\X)\cong \bigoplus^s_{i=0}L^iH^{2i}(\X)$$

\begin{definition}
For $[\alpha]\in \pi_0\widetilde{\V ec}t^s(\X)$,
$\textbf{c}([\alpha])\in \bigoplus^s_{i=0}L^iH^{2i}(\X)$ is called
the total Chern class of $[\alpha]$.
\end{definition}

The inclusions
$$\mathscr{C}^{s, 1}(\P^n_S)(\X) \hookrightarrow \mathscr{C}^{s+1,
1}(\P^{n+1}_S)(\X) \hookrightarrow \mathscr{C}^{s+2,
1}(\P^{n+2}_S)(\X) \hookrightarrow \cdots$$ induce inclusions on
$$\mathscr{C}^{s, 1}(\P^{\infty}_S)(\X) \hookrightarrow \mathscr{C}^{s+1,
1}(\P^{\infty}_S)(\X) \hookrightarrow \mathscr{C}^{s+2,
1}(\P^{\infty}_S)(\X) \hookrightarrow \cdots$$ which induce again maps on
$$\widetilde{\V ec}t^s(\X) \rightarrow \widetilde{\V ec}t^{s+1}(\X)
\rightarrow \widetilde{\V ec}t^{s+2}(\X) \rightarrow \cdots$$ Let
$$\widetilde{\V ec}t(\X):=\lim_{s\to \infty}\widetilde{\V
ec}t^s(\X)$$

\begin{definition}
Suppose that $\X$ is a semi-topological projective variety. Let
$$K_n(\X):=\pi_n\widetilde{\V ec}t(\X)$$ This is called the
$n$-th $K$-group of $\X$.
\end{definition}

This construction of Chern classes is a preparation for a proof of Grothendieck-Riemann-Roch for semi-topological projective varieties.

\begin{example}
When $S=S^0$, $\X=X\times S^0$ where $X$ is a smooth projective variety. Then
$$K_n(\X)=K_n^{semi}(X)$$ where $K^n(X)$ is the semi-topological $K$-group of $X$ constructed by Friedlander and Walker
\cite{FW02, FW03}.
\end{example}

\bibliographystyle{amsplain}

\end{document}